\documentclass[12pt,twoside]{article}

\usepackage[british]{babel}

\usepackage{amsmath}

\DeclareMathOperator{\RE}{Re}
\DeclareMathOperator{\IM}{Im}

\usepackage{amsthm}

\theoremstyle{plain}
\newtheorem{theorem}{Theorem}
\newtheorem{lemma}{Lemma}

\usepackage{amssymb}

\begin{document}

\title{An optimal series expansion of the multiparameter fractional Brownian motion\thanks{This work is supported in part by the Foundation for Knowledge and Competence Development and Sparbanksstiftelsen Nya}}

\author{Anatoliy Malyarenko}

\date{\today}

\maketitle

\begin{abstract}
We derive a series expansion for the multiparameter fractional Brownian motion. The derived expansion is proven to be rate optimal. \textbf{Keywords}: fractional Brownian motion; series expansion; Bessel functions.
\end{abstract}

\section{Introduction}

The \emph{fractional Brownian motion} with Hurst parameter
$H\in(0,1)$ is defined as the centred Gaussian process $\xi(t)$ with the autocorrelation function
\[
R(s,t)=\mathsf{E}\xi(s)\xi(t)=\frac{1}{2}(|s|^{2H}+|t|^{2H}-|s-t|^{2H}).
\]
This process was defined by Kolmogorov \cite{Kol} and became a
popular statistical model after the paper by Mandelbrot and van Ness \cite{Man}.

There exist two multiparameter extensions of the fractional Brownian motion. Both extensions are centred Gaussian random fields on the space $\mathbb{R}^N$. The \emph{multiparameter fractional Brownian sheet} has the autocorrelation function
\[
R(\mathbf{x},\mathbf{y})=\frac{1}{2^N}\prod^N_{j=1}(|x_j|^{2H_j}+|y_j|^{2H_j}-|x_j-y_j|^{2H_j}),\qquad
H_j\in(0,1)
\]
while the \emph{multiparameter fractional Brownian motion} has the autocorrelation function
\begin{equation}\label{mfbmcov}
R(\mathbf{x},\mathbf{y})=\frac{1}{2}(\|\mathbf{x}\|^{2H}+\|\mathbf{y}\|^{2H}
-\|\mathbf{x}-\mathbf{y}\|^{2H}),
\end{equation}
where $\|\cdot\|$ denotes the Euclidean norm in $\mathbb{R}^N$ and where $H\in(0,1)$.

Dzhaparidze and van Zanten \cite{DvZ1} and Igl\'{o}i \cite{Igl}  derived two different explicit series expansions of the fractional Brownian motion. In \cite{DvZ2}, Dzhaparidze and van Zanten extended their previous result to the case of the multiparameter fractional Brownian sheet. All the above mentioned expansions were proven to be rate optimal. We extend the results by Dzhaparidze and van Zanten to the case of the multiparameter fractional Brownian motion.

To present and prove our result, we need to recall definitions of some special functions \cite{Erd}. The \emph{Bessel function of the first kind of order} $\nu$ is defined by the following series
\begin{equation}\label{bessel}
J_{\nu}(z)=\sum^{\infty}_{m=0}\frac{(-1)^mz^{2m+\nu}}{2^mm!\Gamma(\nu+m+1)}.
\end{equation}
Let $j_{\nu,1}<j_{\nu,2}<\dots<j_{\nu,n}<\dots$ be the positive zeros of $J_{\nu}(z)$. Let
\begin{equation}\label{gm}
g_m(u)=2^{(N-2)/2}\Gamma(N/2)\frac{J_{m+(N-2)/2}(u)}{u^{(N-2)/2}},\qquad
m\geq 0,
\end{equation}
and let $\delta^n_m$ denote the Kronecker's delta.

The \emph{Gegenbauer polynomials} of order $m$, $C^{\lambda}_m$, are given by the generating function
\[
\frac{1}{(1-2xt+t^2)^{\lambda}}=\sum^{\infty}_{m=0}C^{\lambda}_m(x)t^m.
\]

Let $m$ be a nonnegative integer, and let $m_0$, $m_1$, \dots, $m_{N-2}$ be integers satisfying the following condition
\[
m=m_0\geq m_1\geq\dots\geq m_{N-2}\geq 0.
\]
Let $\mathbf{x}=(x_1, x_2,\dots, x_N)$ be a point in the space $\mathbb{R}^N$. Let
\[
r_k=\sqrt{x^2_{k+1}+x^2_{k+2}+\dots+x^2_N},
\]
where $k=0$, $1$, \dots, $N-2$. Consider the following functions
\[
\begin{aligned}
H(m_k,\pm,\mathbf{x})&=\left(\frac{x_{N-1}+ix_N}{r_{N-2}}\right)^{\pm m_{N-2}}r^{m_{N-2}}_{N-2}\prod^{N-3}_{k=0}r^{m_k-m_{k+1}}_k\\
&\quad\times C^{m_{k+1}+(N-k-2)/2}_{m_k-m_{k+1}}\left(\frac{x_{k+1}}{r_k}\right),\\
\end{aligned}
\]
and denote
\[
Y(m_k,\pm,\mathbf{x})=r^{-m}_0H(m_k,\pm,\mathbf{x}).
\]
The functions $Y(m_k,\pm,\mathbf{x})$ are called the (complex-valued) \emph{spherical harmonics}. For a fixed $m$, there exist
\begin{equation}\label{hmN}
h(m,N)=\frac{(2m+N-2)(m+N-3)!}{(N-2)!m!}
\end{equation}
spherical harmonics. They are orthogonal in the Hilbert space $L^2(S^{N-1})$ of the square integrable functions on the unit sphere $S^{N-1}$, and the square of the length of the vector $Y(m_k,\pm,\mathbf{x})$ is
\[
L(m_k)=2\pi\prod^{N-2}_{k=1}\frac{\pi2^{k-2m_k-N+2}\Gamma(m_{k-1}+m_k+N-1-k)}
{(m_{k-1}+(N-1-k)/2)(m_{k-1}-m_k)![\Gamma(m_k+(N-1-k)/2)]^2}.
\]

Let $l=l(m_k,\pm)$ be the number of the symbol $(m_0, m_1, \dots, m_{N-2},\pm)$ in the lexicographic ordering. The \emph{real-valued spherical harmonics}, $S^l_m(\mathbf{x})$, can be defined as
\[
S^l_m(\mathbf{x})=
\begin{cases}
Y(m_k,+,\mathbf{x})/\sqrt{L(m_k)},&m_{N-2}=0,\\
\sqrt{2}\RE Y(m_k,+,\mathbf{x})/\sqrt{L(m_k)},&m_{N-2}>0,l=l(m_k,+),\\
-\sqrt{2}\IM Y(m_k,-,\mathbf{x})/\sqrt{L(m_k)},&m_{N-2}>0,l=l(m_k,-).
\end{cases}
\]

The \emph{hypergeometric function} is defined by the series
\[
{}_2F_1(a,b;c;z)=\sum^{\infty}_{k=0}\frac{(a)_k(b)_k}{(c)_kk!}z^k,
\]
where $(u)_k=u(u+1)\dots (u+k-1)$, $(u)_0=u$.

The \emph{incomplete beta function} is defined as
\[
B_z(\alpha,\beta)=\int^z_0t^{\alpha-1}(1-t)^{\beta-1}\,dt.
\]

\begin{theorem}\label{th1}
The multiparameter fractional Brownian motion $\xi(\mathbf{x})$ has the form
\begin{equation}\label{main2}
\xi(\mathbf{x})=\sum^{\infty}_{m=0}\sum^{\infty}_{n=1}\sum^{h(m,N)}_{l=1}\tau_{mn}
[g_m(j_{|m-1|-H,n}\|\mathbf{x}\|)-\delta^0_m]
S^l_m\left(\frac{\mathbf{x}}{\|\mathbf{x}\|}\right)\xi^l_{mn},
\end{equation}
where $\xi^l_{mn}$ are independent standard normal random variables, and
\begin{equation}\label{sigma2}
\tau_{mn}=\frac{2^{H+1}\sqrt{\pi^{(N-2)/2}\Gamma(H+N/2)\Gamma(H+1)\sin(\pi
H)}}{\Gamma(N/2)J_{|m-1|-H+1}(j_{|m-1|-H,n})j^{H+1}_{|m-1|-H,n}}.
\end{equation}
The series \eqref{main2} converges with probability $1$ in the space $C(\mathcal{B})$ of continuous functions in
$\mathcal{B}=\{\,\mathbf{x}\in\mathbb{R}^N\colon\|\mathbf{x}\|\leq
1\,\}$.
\end{theorem}

Ayache and Linde \cite{Aya} derived another representation of the multiparameter fractional Brownian motion, using wavelets. Hence it is natural to compare their result with Theorem~\ref{th1}.

Consider the multiparameter fractional Brownian motion
$\xi(\mathbf{x})$, $\mathbf{x}\in\mathcal{B}$ as the centred
Gaussian random variable $\xi$ on the Banach space $C(\mathcal{B})$ of all continuous functions on $\mathcal{B}$. The $p$th $l$-approximation number of $\xi$ \cite{Li1} is defined by:
\[
l_p(\xi)=\inf\left\{\left(\mathsf{E}\left\|\sum^{\infty}_{j=p}
f_j\xi_j\right\|_{C(\mathcal{B})}^2\right)^{1/2}\colon\xi=\sum^{\infty}_{j=1}f_j\xi_j,
\quad f_j\in C(\mathcal{B})\right\},
\]
where $\xi_j$ are independent standard normal random variables and the infimum is taken over all possible series representations for $\xi$.

Ayache and Linde \cite{Aya} determined the convergence rate of $l_p(\xi)\to 0$ as $p\to\infty$. To formulate their result, introduce the following notation. If $a_n$, $n\geq 1$, and $b_n$, $n\geq 1$ are sequences of positive real numbers, we write $a_n\preceq b_n$ provided that $a_n\leq cb_n$ for a certain $c>0$ and for any positive integer $n$. Then $a_n\approx b_n$ means that $a_n\preceq b_n$ as well as $b_n\preceq a_n$. In the same way, we write $f(\mathbf{u})\preceq g(\mathbf{u})$ provided that $f(\mathbf{u})\leq cg(\mathbf{u})$ for a certain $c>0$ and uniformly for all $\mathbf{u}$, and $f(\mathbf{u})\approx g(\mathbf{u})$ if $f(\mathbf{u})\preceq g(\mathbf{u})$ as well as $g(\mathbf{u})\preceq f(\mathbf{u})$.

Ayache and Linde \cite[Theorem~1.1]{Aya} proved that $l_p(\xi)\approx p^{-H/N}(\log p)^{1/2}$. They also proved that their wavelet series representation possesses the optimal approximation rate. We will prove that our representation \eqref{main2} is optimal as well.

\begin{theorem}\label{th2}
The representation \eqref{main2} possesses the optimal approximation rate for $\xi$ on $\mathcal{B}$.
\end{theorem}

Theorem~\ref{th1} is proved in Section~\ref{proof} while Theorem~\ref{th2} is proved in Section~\ref{rate}.

I am grateful to Professors K.~Dzhaparidze, M.~Lifshits, and H.~van Zanten for useful discussions. I thank the two anonymous referees for their very careful reading of the manuscript and for their helpful remarks.

\section{Proof of Theorem~\ref{th1}}\label{proof}

It is well known \cite{Aya}, that the multiparameter fractional Brownian motion can be represented as the stochastic integral
\[
\xi(\mathbf{x})=c_{HN}\int_{\mathbb{R}^N}\frac{e^{i(\mathbf{p},\mathbf{x})}-1}
{\|\mathbf{p}\|^{N/2+H}}\,d\hat{W}(\mathbf{p}),
\]
where $d\hat{W}$ is the complex-valued white noise obtained by Fourier transformation of the real-valued white noise.  It follows that the covariance function \eqref{mfbmcov} of the multiparameter fractional Brownian motion can be represented as
\begin{equation}\label{mfbmcov2}
R(\mathbf{x},\mathbf{y})=c^2_{HN}\int_{\mathbf{R}^N}[e^{i(\mathbf{p},\mathbf{x})}-1]
[e^{-i(\mathbf{p},\mathbf{y})}-1]\|\mathbf{p}\|^{-N-2H}\,d\mathbf{p}.
\end{equation}
The following Lemma gives the explicit value of the constant $c^2_{HN}$. This result was announced by Malyarenko \cite{Mal}.

\begin{lemma}
The constant $c^2_{HN}$ has the following value:
\[
c^2_{HN}=\frac{2^{2H-1}\Gamma(H+N/2)\Gamma(H+1)\sin(\pi
H)}{\pi^{(N+2)/2}}.
\]
\end{lemma}

\begin{proof}
For $N=1$, our formula has the form
\[
c^2_{H1}=\frac{2^{2H-1}\Gamma(H+1/2)\Gamma(H+1)\sin(\pi
H)}{\pi^{3/2}},
\]
or
\[
c^2_{H1}=\frac{\Gamma(2H+1)\sin(\pi H)}{2\pi}.
\]
Here we used the doubling formula for gamma function. This result is known \cite{ST}. Therefore, in the rest of the proof we can and will suppose that $N\geq 2$.

Rewrite \eqref{mfbmcov2} as
\begin{equation}\label{mfbmcov3}
\begin{aligned}
R(\mathbf{x},\mathbf{y})&=c^2_{HN}\int_{\mathbf{R}^N}[1-e^{i(\mathbf{p},\mathbf{x})}]
\|\mathbf{p}\|^{-N-2H}\,d\mathbf{p}\\
&\quad+c^2_{HN}\int_{\mathbf{R}^N}[1-e^{-i(\mathbf{p},\mathbf{y})}]
\|\mathbf{p}\|^{-N-2H}\,d\mathbf{p}\\
&\quad-c^2_{HN}\int_{\mathbf{R}^N}[1-e^{i(\mathbf{p},\mathbf{x}-\mathbf{y})}]
\|\mathbf{p}\|^{-N-2H}\,d\mathbf{p}.
\end{aligned}
\end{equation}
Consider the first term in the right hand side of this formula.
Using formula 3.3.2.3 from \cite{Pru}, we obtain
\begin{multline*}
c^2_{HN}\int_{\mathbf{R}^N}[1-e^{i(\mathbf{p},\mathbf{x})}]
\|\mathbf{p}\|^{-N-2H}\,d\mathbf{p}\\
=\frac{2\pi^{(N-1)/2}c^2_{HN}}{\Gamma((N-1)/2)}\int^{\infty}_0\lambda^{N-1}\,
d\lambda\int^{\pi}_0[1-e^{i\lambda\|\mathbf{x}\|\cos
u}]\lambda^{-N-2H}\sin^{N-2}u\,du.
\end{multline*}

It is clear that the integral of the imaginary part is equal to $0$.
The integral of the real part may be rewritten as
\begin{multline}\label{realint}
c^2_{HN}\int_{\mathbf{R}^N}[1-e^{i(\mathbf{p},\mathbf{x})}]
\|\mathbf{p}\|^{-N-2H}\,d\mathbf{p}\\
=\frac{2\pi^{(N-1)/2}c^2_{HN}}{\Gamma((N-1)/2)}\int^{\infty}_0\lambda^{-1-2H}
\int^{\pi}_0[1-\cos(\lambda\|\mathbf{x}\|\cos
u)]\sin^{N-2}u\,du\,d\lambda.
\end{multline}

To calculate the inner integral, we use formulas 2.5.3.1 and
2.5.55.7 from \cite{Pru}:
\[
\begin{aligned}
\int^{\pi}_0\sin^{N-2}u\,du&=\frac{\sqrt{\pi}\Gamma((N-1)/2)}{\Gamma(N/2)},\\
\int^{\pi}_0\cos(\lambda\|\mathbf{x}\|\cos
u)\sin^{N-2}u\,du&=\sqrt{\pi}2^{(N-2)/2}\Gamma((N-1)/2)
\frac{J_{(N-2)/2}(\lambda\|\mathbf{x}\|)}
{(\lambda\|\mathbf{x}\|)^{(N-2)/2}}.
\end{aligned}
\]
It follows that
\begin{equation}\label{innerint}
\int^{\pi}_0[1-\cos(\lambda\|\mathbf{x}\|\cos u)]\sin^{N-2}u\,du
=\frac{\sqrt{\pi}\Gamma((N-1)/2)}{\Gamma(N/2)}[1-g_0(\lambda\|\mathbf{x}\|)].
\end{equation}
Substituting \eqref{innerint} in \eqref{realint}, we obtain
\begin{equation}\label{realint4}
c^2_{HN}\int_{\mathbf{R}^N}[1-e^{i(\mathbf{p},\mathbf{x})}]
\|\mathbf{p}\|^{-N-2H}\,d\mathbf{p}
=\frac{2\pi^{N/2}c^2_{HN}}{\Gamma(N/2)}\int^{\infty}_0\lambda^{-2H-1}
[1-g_0(\lambda\|\mathbf{x}\|)]\,d\lambda.
\end{equation}

To calculate this integral, we use formula 2.2.3.1 from \cite{Pru}:
\[
\int^1_0(1-v^2)^{(N-3)/2}\,dv=\frac{\sqrt{\pi}\Gamma((N-1)/2)}{2\Gamma(N/2)}.
\]
It follows that
\begin{equation}\label{p1}
1=\frac{2\Gamma(N/2)}{\sqrt{\pi}\Gamma((N-1)/2)}\int^1_0(1-v^2)^{(N-3)/2}\,dv.
\end{equation}
On the other hand, according to formula 2.5.6.1 from \cite{Pru} we
have
\[
\int^1_0(1-v^2)^{(N-3)/2}\cos(\lambda\|\mathbf{x}\|v)\,dv=\sqrt{\pi}
2^{(N-4)/2}\Gamma((N-1)/2)\frac{J_{(N-2)/2}(\lambda\|\mathbf{x}\|)}
{(\lambda\|\mathbf{x}\|)^{(N-2)/2}}.
\]
It follows that
\begin{equation}\label{p2}
g_0(\lambda\|\mathbf{x}\|)=\frac{2\Gamma(N/2)}{\sqrt{\pi}\Gamma((N-1)/2)}
\int^1_0(1-v^2)^{(N-3)/2}\cos(\lambda\|\mathbf{x}\|v)\,dv.
\end{equation}
Subtracting \eqref{p2} from \eqref{p1}, we obtain
\[
1-g_0(\lambda\|\mathbf{x}\|)=\frac{4\Gamma(N/2)}{\sqrt{\pi}\Gamma((N-1)/2)}
\int^1_0(1-v^2)^{(N-3)/2}\sin^2\left(\frac{\lambda\|\mathbf{x}\|}{2}v\right)\,dv.
\]
Substitute this formula in \eqref{realint4}. We have
\begin{multline}\label{doubleint}
c^2_{HN}\int_{\mathbf{R}^N}[1-e^{i(\mathbf{p},\mathbf{x})}]
\|\mathbf{p}\|^{-N-2H}\,d\mathbf{p}\\
=\frac{8\pi^{(N-1)/2}c^2_{HN}}{\Gamma((N-1)/2)}\int^{\infty}_0\int^1_0
\lambda^{-2H-1}(1-v^2)^{(N-3)/2}
\sin^2\left(\frac{\lambda\|\mathbf{x}\|}{2}v\right)\,dv\,d\lambda.
\end{multline}

Consider two integrals:
\[
\int^1_0\lambda^{-2H-1}\sin^2\left(\frac{v\|\mathbf{x}\|}{2}\lambda\right)\,d\lambda
\quad\text{and}\quad\int^{\infty}_1\lambda^{-2H-1}\sin^2\left(
\frac{v\|\mathbf{x}\|}{2}\lambda\right)\,d\lambda.
\]
In the first integral, we bound the second multiplier by
$\lambda^2/4$. In the second integral, we bound it by $1$. It
follows that the integral
\[
\int^{\infty}_0\lambda^{-2H-1}\sin^2\left(\frac{v\|\mathbf{x}\|}{2}\lambda\right)\,d\lambda
\]
converges uniformly, and we are allowed to change the order of
integration in the right hand side of \eqref{doubleint}. After that
the inner integral is calculated using formula 2.5.3.13 from
\cite{Pru}:
\[
\int^{\infty}_0\lambda^{-2H-1}\sin^2\left(\frac{v\|\mathbf{x}\|}{2}\lambda\right)
\,d\lambda=\frac{\pi v^{2H}}{4\Gamma(2H+1)\sin(\pi
H)}\|\mathbf{x}\|^{2H}.
\]
Now formula \eqref{realint4} may be rewritten as
\begin{multline*}
c^2_{HN}\int_{\mathbf{R}^N}[1-e^{i(\mathbf{p},\mathbf{x})}]
\|\mathbf{p}\|^{-N-2H}\,d\mathbf{p}\\
=\frac{2\pi^{(N-1)/2}c^2_{HN}}{\Gamma((N-1)/2)\Gamma(2H+1)\sin(\pi
H)}\int^1_0v^{2H}(1-v^2)^{(N-3)/2}\,dv\cdot\|\mathbf{x}\|^{2H}.
\end{multline*}
The integral in the right hand side can be calculated by formula
2.2.4.8 from \cite{Pru}:
\[
\int^1_0v^{2H}(1-v^2)^{(N-3)/2}\,dv=\frac{\Gamma((N-1)/2)\Gamma(H+1/2)}{2\Gamma(H+N/2)}
\]
and \eqref{realint4} is rewritten once more as
\begin{multline}\label{realint5}
c^2_{HN}\int_{\mathbf{R}^N}[1-e^{i(\mathbf{p},\mathbf{x})}]
\|\mathbf{p}\|^{-N-2H}\,d\mathbf{p}\\
=\frac{\pi^{(N+2)/2}c^2_{HN}}{2^{2H}\Gamma(H+1)\Gamma(H+N/2)\sin(\pi
H)}\|\mathbf{x}\|^{2H}.
\end{multline}

On the other hand, the left hand side of \eqref{realint5} is clearly
equal to $\dfrac{1}{2}\|\mathbf{x}\|^{2H}$. The statement of the Lemma follows.
\end{proof}

Taking into account \eqref{realint4}, one can rewrite
\eqref{mfbmcov3} as
\begin{equation}\label{mfbmcov4}
R(\mathbf{x},\mathbf{y})=\frac{2\pi^{N/2}c^2_{HN}}{\Gamma(N/2)}
\int^{\infty}_0\lambda^{-1-2H}[1-g_0(\lambda\|\mathbf{x}\|)
-g_0(\lambda\|\mathbf{y}\|)+g_0(\lambda\|\mathbf{x}-\mathbf{y}\|)]\,d\lambda.
\end{equation}
Let $\mathbf{x}\neq\mathbf{0}$, $\mathbf{y}\neq\mathbf{0}$ and let $\varphi$ denote the angle between the vectors $\mathbf{x}$ and $\mathbf{y}$. Two addition theorems for Bessel functions (formulas 7.15(30) and 7.15(31) from \cite[vol.~2]{Erd}) may be written in our notation as
\begin{equation}\label{addition}
g_0(\lambda\|\mathbf{x}-\mathbf{y}\|)=
\sum^{\infty}_{m=0}h(m,N)g_m(\lambda\|\mathbf{x}\|)
g_m(\lambda\|\mathbf{y}\|)\frac{C^{(N-2)/2}_m(\cos\varphi)}{C^{(N-2)/2}_m(1)},
\end{equation}
where $g_m$ is as in \eqref{gm}. Substituting \eqref{addition} in \eqref{mfbmcov4}, we obtain
\begin{equation}\label{covkappa}
R(\mathbf{x},\mathbf{y})=\frac{2\pi^{N/2}}{\Gamma(N/2)}\sum^{\infty}_{m=0}
\frac{C^{(N-2)/2}_m(\cos\varphi)}{C^{(N-2)/2}_m(1)}h(m,N)
\int^{\infty}_0\kappa^m_{\|\mathbf{x}\|}(\lambda)\kappa^m_{\|\mathbf{y}\|}(\lambda)\,
d\lambda,
\end{equation}
where
\[
\kappa^m_r(\lambda)=c_{HN}\lambda^{-H-1/2}[g_m(r\lambda)-\delta^0_m].
\]

Recall that the \emph{Hankel transform} of order $\nu>-1$ of a function $\kappa\in L_2(0,\infty)$ (see \cite[Section~8.4]{Tit}) is defined as
\[
\hat{\kappa}(u)=\int^{\infty}_0\kappa(\lambda)J_{\nu}(u\lambda)\sqrt{u\lambda}\,d\lambda.
\]
Define the kernel $\hat{\kappa}^m_r(u)$ as the Hankel transform of order $|m-1|-H$ of the kernel $\kappa^m_r(\lambda)$:
\begin{equation}\label{ht}
\hat{\kappa}^m_r(u)=\int^{\infty}_0\kappa^m_r(\lambda)J_{|m-1|-H}(u\lambda)\sqrt{u\lambda}\,d\lambda.
\end{equation}
Therefore, the Parceval identity for the Hankel transform \cite[Section~8.5, Theorem~129]{Tit} implies that \eqref{covkappa} may be rewritten as
\begin{equation}\label{covk}
R(\mathbf{x},\mathbf{y})=\frac{2\pi^{N/2}}{\Gamma(N/2)}\sum^{\infty}_{m=0}
\frac{C^{(N-2)/2}_m(\cos\varphi)}{C^{(N-2)/2}_m(1)}h(m,N)
\int^{\infty}_0\hat{\kappa}^m_{\|\mathbf{x}\|}(u)\hat{\kappa}^m_{\|\mathbf{y}\|}(u)\,du.
\end{equation}

To calculate $\hat{\kappa}^m_r(u)$ in the case of $m\geq 1$, we use formula 2.12.31.1 from \cite{Pru2}.
\begin{equation}\label{km}
\hat{\kappa}^m_r(u)=\frac{c_{HN}u^{-H+m-1/2}(r^2-u^2)^{H+N/2-1}}{2^{H+N/2-1}\Gamma(H+N/2)r^{m+N-2}}
\chi_{(0,r)}(u),
\end{equation}
where $\chi_{(0,r)}(u)$ denote the indicator function of the
interval $(0,r)$.

For $m=0$, the integral in \eqref{ht} can be rewritten as
\[
\begin{aligned}
\hat{\kappa}^0_r(u)&=2^{(N-2)/2}\Gamma(N/2)c_{HN}r^{1-N/2}\sqrt{u}\int^{\infty}_0\lambda^{1-H-N/2}
J_{(N-2)/2}(r\lambda)J_{1-H}(u\lambda)\,d\lambda\\
&\quad-\sqrt{u}c_{HN}\int^{\infty}_0\lambda^{-H}J_{1-H}(u\lambda)\,d\lambda.
\end{aligned}
\]
To calculate the second integral, we use formula~2.12.2.2 from
\cite{Pru2}.
\[
\int^{\infty}_0\lambda^{-H}J_{1-H}(u\lambda)\,d\lambda=\frac{\Gamma(1-H)}{2^Hu^{1-H}}.
\]
For the first integral, we use formula 2.12.31.1 from \cite{Pru2} once more. In the case of $u>r$ we obtain
\[
\int^{\infty}_0\lambda^{1-H-N/2}J_{(N-2)/2}(r\lambda)J_{1-H}(u\lambda)\,d\lambda
=\frac{\Gamma(1-H)r^{(N-2)/2}}{2^{H+(N-2)/2}\Gamma(N/2)u^{1-H}}.
\]
In the case of $u<r$, we have
\[
\begin{aligned}
\int^{\infty}_0\lambda^{1-H-N/2}J_{(N-2)/2}(r\lambda)J_{1-H}(u\lambda)\,d\lambda
&=\frac{\Gamma(1-H)r^{2H+N/2-3}}{2^{H+(N-2)/2}\Gamma(H+(N-2)/2)\Gamma(2-H)}\\
&\quad\times{}_2F_1(1-H,2-H-N/2;2-H;u^2/r^2).
\end{aligned}
\]
Using formula 7.3.1.28 from \cite{Pru3}, we can rewrite the last expression as
\[
\begin{aligned}
\int^{\infty}_0\lambda^{1-H-N/2}J_{(N-2)/2}(r\lambda)J_{1-H}(u\lambda)\,d\lambda
&=\frac{r^{(N-2)/2}}{2^{H+(N-2)/2}\Gamma(H+(N-2)/2)u^{1-H}}\\
&\quad\times B_{u^2/r^2}(1-H,H+(N-2)/2).
\end{aligned}
\]

Combining everything together, we obtain
\[
\begin{aligned}
\hat{\kappa}^0_r(u)&=2^{-H}u^{H-1/2}c_{HN}[-\Gamma(1-H)+\frac{\Gamma(N/2)}{\Gamma(H+(N-2)/2)}\\
&\quad\times B_{u^2/r^2}(1-H,H+(N-2)/2)]\chi_{(0,r)}(u).
\end{aligned}
\]
It follows from the last display and from \eqref{km} that the support of the kernels $\hat{\kappa}^m_{\|\mathbf{x}\|}(u)$ and $\hat{\kappa}^m_{\|\mathbf{y}\|}(u)$ lies in $[0,1]$ since $\|\mathbf{x}\|$, $\|\mathbf{y}\|\leq 1$ for $\mathbf{x}$, $\mathbf{y}\in\mathcal{B}$. We rewrite \eqref{covk} as
\begin{equation}\label{covb}
R(\mathbf{x},\mathbf{y})=\frac{2\pi^{N/2}}{\Gamma(N/2)}\sum^{\infty}_{m=0}
\frac{C^{(N-2)/2}_m(\cos\varphi)}{C^{(N-2)/2}_m(1)}h(m,N)
\int^1_0k^m_{\|\mathbf{x}\|}(u)k^m_{\|\mathbf{y}\|}(u)\,du
\end{equation}
for $\mathbf{x}$, $\mathbf{y}\in\mathcal{B}$.

For any $\nu>-1$, the Fourier--Bessel functions
\[
\varphi_{\nu,n}(u)=\frac{\sqrt{2u}}{J_{\nu+1}(j_{\nu,n})}J_{\nu}(j_{\nu,n}u),\qquad
n\geq 1
\]
form a complete, orthonormal system in $L^2[0,1]$ \cite[Section~18.24]{Wat}. Put $\nu=|m-1|-H$. We obtain
\[
\int^1_0\hat{\kappa}^m_{\|\mathbf{x}\|}(u)\hat{\kappa}^m_{\|\mathbf{y}\|}(u)\,du=
\sum^{\infty}_{n=1}b^m_n(\|\mathbf{x}\|)b^m_n(\|\mathbf{y}\|),
\]
where
\[
b^m_n(\|\mathbf{x}\|)=\int^{1}_0\hat{\kappa}^m_{\|\mathbf{x}\|}(u)\varphi_{|m-1|-H,n}(u)\,du.
\]
This integral can be rewritten as
\[
b^m_n(\|\mathbf{x}\|)=\int^{\infty}_0\hat{\kappa}^m_{\|\mathbf{x}\|}(u)\varphi_{|m-1|-H,n}(u)\,du.
\]
To calculate it, we use the inversion formula for Hankel transform:
\[
\int^{\infty}_0\hat{\kappa}^m_r(u)J_{\nu}(\lambda u)\sqrt{\lambda u}\,du=\kappa^m_r(\lambda)
\]
and obtain
\[
b^m_n(\|\mathbf{x}\|)=\frac{\sqrt{2}}{J_{|m-1|-H+1}(j_{|m-1|-H,n})\sqrt{j_{|m-1|-H,n}}}
\kappa^m_{\|\mathbf{x}\|}(j_{|m-1|-H,n}).
\]
Substituting the calculated values in \eqref{covb}, we obtain
\begin{equation}\label{Rxy}
\begin{aligned}
R(\mathbf{x},\mathbf{y})&=\frac{4\pi^{N/2}}{\Gamma(N/2)}\sum^{\infty}_{m=0}\sum^{\infty}_{n=1}
\frac{C^{(N-2)/2}_m(\cos\varphi)}{C^{(N-2)/2}_m(1)}h(m,N)\\
&\quad\times\frac{\kappa^m_{\|\mathbf{x}\|}(j_{|m-1|-H,n})
\kappa^m_{\|\mathbf{y}\|}(j_{|m-1|-H,n})}{J^2_{|m-1|-H+1}(j_{|m-1|-H,n})j_{|m-1|-H,n}}.
\end{aligned}
\end{equation}

According to the addition theorem for spherical harmonics
\cite[Vol.~2, Chapter~XI, section~4, Theorem~4]{Erd},
\begin{equation}\label{add}
\frac{C^{(N-2)/2}_m(\cos\varphi)}{C^{(N-2)/2}_m(1)}=\frac{2\pi^{N/2}}
{\Gamma(N/2)h(m,N)}\sum^{h(m,N)}_{l=1}S^l_m\left(\frac{\mathbf{x}}{\|\mathbf{x}\|}\right)
S^l_m\left(\frac{\mathbf{y}}{\|\mathbf{y}\|}\right).
\end{equation}
Substituting this equality in \eqref{Rxy}, we obtain
\[
\begin{aligned}
R(\mathbf{x},\mathbf{y})&=\frac{8\pi^N}{(\Gamma(N/2))^2}
\sum^{\infty}_{m=0}\sum^{\infty}_{n=1}\sum^{h(m,N)}_{l=1}
\frac{\kappa^m_{\|\mathbf{x}\|}(j_{|m-1|-H,n})
\kappa^m_{\|\mathbf{y}\|}(j_{|m-1|-H,n})}
{J^2_{|m-1|-H+1}(j_{|m-1|-H,n})j_{|m-1|-H,n}}\\
&\quad\times S^l_m\left(\frac{\mathbf{x}}{\|\mathbf{x}\|}\right)
S^l_m\left(\frac{\mathbf{y}}{\|\mathbf{y}\|}\right).
\end{aligned}
\]
It follows immediately that the random field $\xi(\mathbf{x})$
itself has the form \eqref{main2}, \eqref{sigma2}, and the series
\eqref{main2} converges in mean square for any fixed
$\mathbf{x}\in\mathcal{B}$.

Since the functions in \eqref{main2} are continuous and the random variables are symmetric and independent, the It\^{o}--Nisio theorem \cite{Vak} implies that for proving that the series \eqref{main2} converges uniformly on $\mathcal{B}$ with probability $1$, it is sufficient to show that the corresponding sequence of partial sums
\begin{equation}\label{partialsums}
\begin{aligned}
\xi_M(\mathbf{x})=\sum^M_{m=0}\sum^M_{n=1}\sum^{h(m,N)}_{l=1}\tau_{mn}
[g_m(j_{|m-1|-H,n}\|\mathbf{x}\|)-\delta^0_m]
S^l_m\left(\frac{\mathbf{x}}{\|\mathbf{x}\|}\right)\xi^l_{mn}
\end{aligned}
\end{equation}
is weakly relatively compact in the space $C(\mathcal{B})$. To prove this, one can use the same method as in the proof of Theorem~4.1 in \cite{DvZ3}. Proof of Theorem~\ref{th1} is finished.

Consider some particular cases of Theorem~\ref{th1}. In the case of $N=1$, there exists only two spherical harmonics on the $0$-dimensional sphere $S^0=\{-1,1\}$. They are
\[
S^1_0(x)=\frac{1}{\sqrt{2}},\qquad S^1_1(x)=\frac{x}{\sqrt{2}},
\]
where $x\in S^0$. Moreover, we have
\[
g_0(j_{1-H,n}|x|)=\cos(j_{1-H,n}|x|),\qquad
g_1(j_{-H,n}|x|)=\sin(j_{-H,n}|x|).
\]
Therefore formula \eqref{main2} becomes
\begin{equation}\label{DzvZa}
\xi(x)=2c_{H1}\left(\sum^{\infty}_{n=1}\frac{\cos(j_{1-H,n}x)-1}
{J_{2-H}(j_{1-H,n})j^{1+H}_{1-H,n}}\xi^1_{0n}+\sum^{\infty}_{n=1}\frac{\sin(j_{-H,n}x)}
{J_{1-H}(j_{-H,n})j^{1+H}_{-H,n}}\xi^1_{1n}\right).
\end{equation}
This is the result of \cite{DvZ1}.

Consider the case of the multiparameter fractional Brownian motion on the plane ($N=2$). Let $(r,\varphi)$ be the polar coordinates. The spherical harmonics are:
\[
S^1_0(\varphi)=\frac{1}{\sqrt{2\pi}},\quad
S^1_m(\varphi)=\frac{\cos(m\varphi)}{\sqrt{\pi}},\quad
S^2_m(\varphi)=\frac{\sin(m\varphi)}{\sqrt{\pi}}.
\]
It follows that
\[
\begin{aligned}
\xi(r,\varphi)&=\frac{2^{H+1}\Gamma(H+1)\sqrt{\sin(\pi
H)}}{\sqrt{\pi}}\left(\frac{1}{\sqrt{2}}\sum^{\infty}_{n=1}
\frac{J_0(j_{1-H,n}r)-1}{J_{2-H}(j_{1-H,n})j^{1+H}_{1-H,n}}\xi^1_{0n}\right.\\
&\quad+\sum^{\infty}_{m=1}\sum^{\infty}_{n=1}
\frac{J_m(j_{m-1-H,n}r)\cos(m\varphi)}{J_{m-H}(j_{m-1-H,n})j^{1+H}_{m-1-H,n}}\xi^1_{mn}\\
&\left.\quad+\sum^{\infty}_{m=1}\sum^{\infty}_{n=1}
\frac{J_m(j_{m-1-H,n}r)\sin(m\varphi)}{J_{m-H}(j_{m-1-H,n})j^{1+H}_{m-1-H,n}}\xi^2_{mn}
\right).
\end{aligned}
\]

\section{Proof of Theorem~\ref{th2}}\label{rate}

It is enough to prove that the rate of convergence in \eqref{main2} is not more than the optimal rate $p^{-H/N}(\log p)^{1/2}$, where $p$ denote the number of terms in a suitable truncation of the series. Since for $N>1$ we have a triple sum in our expansion \eqref{main2}, it is not clear a priori how we should truncate the series. We need a lemma.

\begin{lemma}\label{l2}
The Bessel functions, the functions $g_m$ defined by \eqref{gm}, the numbers $h(m,N)$, and the spherical harmonics $S^m_l(\mathbf{x}/\|\mathbf{x}\|)$ have the following properties.
\begin{enumerate}
\item\label{i1} $j_{\nu,n}\approx n+\nu/2-1/4$.
\item\label{i2} $J^2_{\nu+1}(j_{\nu,n})\approx\frac{1}{j_{\nu,n}}$.
\item\label{i3} $|g_0(u)|\leq 1$.
\item\label{i4} $|g_m(u)|\preceq\frac{1}{[m(m+N-2)]^{(N-1)/4}}$, $m\geq 1$, $N\geq 2$.
\item\label{i5} $|g'_m(u)|\preceq\frac{1}{[m(m+N-2)]^{(N-3)/4}}$, $m\geq 1$, $N\geq 2$.

\item\label{i6} $h(m,N)\preceq m^{N-2}$, $n\geq 2$.

\item\label{i7} $|S^l_m(\mathbf{x}/\|\mathbf{x}\|)|\preceq m^{(N-2)/2}$, $N\geq 2$.
\end{enumerate}
\end{lemma}

\begin{proof}
Property~\ref{i1} is proved in \cite[Section~15.53]{Wat}. It is shown in \cite[Section~7.21]{Wat} that
\begin{equation}\label{Wat}
J^2_{\nu}(u)+J^2_{\nu+1}(u)\approx\frac{1}{u},
\end{equation}
since Property~\ref{i2}. Property~\ref{i3} follows from the fact that $g_0(u)$ is the element of the unitary matrix \cite{Vil}.

Let $\mu_1<\mu_2<\dots$ be the sequence of positive stationary values of the Bessel function $J_{m+(N-2)/2}(u)$. It is shown in \cite[Section~15.31]{Wat} that $|J_{m+(N-2)/2}(\mu_1)|>|J_{m+(N-2)/2}(\mu_2)|>\dots$.

Let $u_1$ be the first positive maximum of the function $g_m(u)$. Let $x_1$ denote the maximal value of the function $|g_m(u)|$ in the interval $(0,j_{m+(N-2)/2,1})$, let $x_2$ denote the maximal value of $|g_m(u)|$ in the interval $(j_{m+(N-2)/2,1},j_{m+(N-2)/2,2})$, and so on. Then we have
\[
\begin{aligned}
x_1&\preceq\frac{|J_{m+(N-2)/2}(\mu_1)|}{u_1^{(N-2)/2}},\\
x_2&\preceq\frac{|J_{m+(N-2)/2}(\mu_2)|}{j_{m+(N-2)/2,1}^{(N-2)/2}},\\
x_3&\preceq\frac{|J_{m+(N-2)/2}(\mu_3)|}{j_{m+(N-2)/2,2}^{(N-2)/2}},
\end{aligned}
\]
and so on. The right hand sides form a decreasing sequence. Then
\[
|g_m(u)|\preceq\frac{|J_{m+(N-2)/2}(\mu_1)|}{u_1^{(N-2)/2}}.
\]
It follows from \eqref{Wat} that
\[
|J_{\nu}(u)|\preceq\frac{1}{\sqrt{u}},
\]
and we obtain
\[
|g_m(u)|\preceq\frac{1}{\mu_1u_1^{(N-2)/2}}
\]

To estimate $u_1$, consider the differential equation
\begin{equation}\label{diff}
u^2\frac{d^2f}{du^2}+(N-1)u\frac{df}{du}+[u^2-m(m+N-2)]f=0.
\end{equation}
It follows from formulas (3) and (4) in \cite[Section~4.31]{Wat}
that the function $g_m$ satisfies this equation. It follows from the series representation \eqref{bessel} of the Bessel function and the definition \eqref{gm} of the function $g_m(u)$ that for $m\geq 1$ we have $g_m(0)=0$ and $g_m$ increases in some right neighbourhood of zero. Therefore we have $g_m(u_1)>0$, $g'_m(u_1)=0$ and $g''_m(u_1)\leq 0$. It follows from equation \eqref{diff} that $u_1\geq\sqrt{m(m+N-2)}$.

For $\mu_1$, we have the estimate $\mu_1>m+(N-2)/2$ (\cite[Section~15.3]{Wat}). Using the inequality $m+(N-2)/2\geq\sqrt{m(m+N-2)}$, we obtain Property~\ref{i4}.

In any extremum of $g'_m(u)$ we have $g''_m(u)=0$. It follows from
\eqref{diff} that
\[
|g'_m(u)|=\left|\frac{[u^2-m(m+N-2)]g_m(u)}{(N-1)u}\right|\preceq|ug_m(u)|.
\]
Property~\ref{i5} follows from this formula and Property~\ref{i4}.

Property~\ref{i6} follows from \eqref{hmN} and Stirling's formula.

It follows from \eqref{add} that
\[
\left(S^l_m\left(\frac{\mathbf{x}}{\|\mathbf{x}\|}\right)\right)^2\leq
\frac{\Gamma(N/2)h(m,N)}{2\pi^{N/2}}\frac{C^{(N-2)/2}_m(\cos\varphi)}{C^{(N-2)/2}_m(1)}.
\]
By \cite[Vol.~2, formula~10.18(7)]{Erd},
\[
\max_{0\leq\varphi\leq\pi}\left|C^{(N-2)/2}_m(\cos\varphi)\right|=C^{(N-2)/2}_m(1).
\]
It follows that $|S^l_m(\mathbf{x}/\|\mathbf{x}\|)|\preceq\sqrt{h(m,N)}$. Now Property~\ref{i7} follows from Property~\ref{i6}.
\end{proof}

Denote
\begin{equation}\label{u}
u^l_{mn}(\mathbf{x})=\tau_{mn}[g_m(j_{|m-1|-H,n}\|\mathbf{x}\|)-\delta^0_m]
S^l_m\left(\frac{\mathbf{x}}{\|\mathbf{x}\|}\right).
\end{equation}
Using Lemma~\ref{l2}, we can write the following estimate:
\begin{equation}\label{est}
\mathsf{E}\left[\sum^{h(m,N)}_{l=1}u^l_{mn}(\mathbf{x})\xi^l_{mn}\right]^2
\preceq\frac{1}{(m+1)(m/2+n)^{2H+1}}.
\end{equation}
For any positive integer $q$, consider the truncation
\[
\sum_{(m+1)(m/2+n)^{2H+1}\leq
q}\sum^{h(m,N)}_{l=1}u^l_{mn}(\mathbf{x})\xi^l_{mn}
\]
of the expansion \eqref{main2}. The number of terms in this sum is asymptotically equal to the integral
\[
\iint_{(u+1)(u/2+v)\leq q,u\geq 0,v\geq 1}u^{N-2}\,du\,dv
\]
and therefore is bounded by a constant times $q^{N/(2H+2)}$. Therefore we need to prove that
\[
\begin{aligned}
&\left(\mathsf{E}\left\|\sum_{(m+1)(m/2+n)^{2H+1}>
q}\sum^{h(m,N)}_{l=1}u^l_{mn}(\mathbf{x})\xi^l_{mn}\right\|_{C(\mathcal{B})}^2\right)^{1/2}\\
&\quad\preceq q^{(N/(2H+2))(-H/N)}(\log q^{N/(2H+2)})^{1/2}
\preceq q^{-H/(2H+2)}(\log q)^{1/2}.
\end{aligned}
\]
By equivalence of moments (\cite[Proposition~2.1]{Khn}), the last formula is equivalent to the following asymptotic relation.
\[
\mathsf{E}\sup_{\mathbf{x}\in\mathcal{B}}\left|\sum_{(m+1)(m/2+n)^{2H+1}>
q}\sum^{h(m,N)}_{l=1}u^l_{mn}(\mathbf{x})\xi^l_{mn}\right|\preceq
q^{-H/(2H+2)}(\log q)^{1/2}.
\]

To prove this relation, consider the partial sum $\eta_k(\mathbf{x})$ defined by
\[
\eta_k(\mathbf{x})=\sum_{2^{k-1}<(m+1)(m/2+n)^{2H+1}\leq
2^k}\sum^{h(m,N)}_{l=1}u^l_{mn}(\mathbf{x})\xi^l_{mn}.
\]
For a given $\varepsilon_k>0$, let $\mathbf{x}_1$, \dots, $\mathbf{x}_{P_k}\in S^{N-1}$ be a maximal $\varepsilon_k$-net in $S^{N-1}$, i.e., the angle $\varphi(\mathbf{x}_j,\mathbf{x}_k)$ between any two different vectors $\mathbf{x}_j$ and $\mathbf{x}_k$ is greater than $\varepsilon_k$ and the addition of any new point breaks this property. Then $P_k\approx\varepsilon_k^{-N+1}$. A proof of this fact for the case of $N=3$ by Baldi et al \cite[Lemma~5]{Bal}, is easy generalised to higher dimensions. The \emph{Voronoi cell} $S(\mathbf{x}_j)$ is defined as
\[
S(\mathbf{x}_j)=\{\,\mathbf{x}\in S^{N-1}\colon\varphi(\mathbf{x},\mathbf{x}_j)\leq\varphi(\mathbf{x},\mathbf{x}_k), k\neq j\,\}.
\]
Divide the ball $\mathcal{B}$ onto $[\varepsilon^{-1}_k]$ concentric spherical layers of thickness $\preceq\varepsilon_k$. Voronoi cells determine the division of each layer onto $P_k$~segments. The angle $\varphi(\mathbf{x},\mathbf{y})$ between any two vectors $\mathbf{x}$ and $\mathbf{y}$ in the same segment is $\preceq\varepsilon_k$ (we never choose the point $\mathbf{0}$). Call all segments in all layers $\mathcal{B}_j$, $1\leq j\leq M_k=P_k[\varepsilon^{-1}_k]\preceq\varepsilon_k^{-N}$. Choose a point $\mathbf{x}_j$ inside each segment. Then we have
\[
\mathsf{E}\sup_{\mathbf{x}\in\mathcal{B}}|\eta_k(\mathbf{x})|\leq
\mathsf{E}\sup_{1\leq j\leq M_k}|\eta_k(\mathbf{x}_j)|+
\mathsf{E}\sup_{1\leq j\leq
M_k}\sup_{\mathbf{x},\mathbf{y}\in\mathcal{B}_j}|\eta_k(\mathbf{x})-\eta_k(\mathbf{y})|.
\]

By a maximal inequality for Gaussian sequences
\cite[Lemma~2.2.2]{Waa}, the first term is bounded by a positive constant times
\[
\sqrt{1+\log M_k}\sup_{1\leq j\leq
M_k}\sqrt{\mathsf{E}(\eta_k(\mathbf{x}_j))^2}.
\]
Using the estimate \eqref{est}, we obtain
\[
\mathsf{E}(\eta_k(\mathbf{x}))^2\preceq\iint_{{2^{k-1}<(u+1)(u/2+v)^{2H+1}\leq
2^k,u\geq 0,v\geq 1}}\frac{du\,dv}{(u+1)(u/2+v)^{2H+1}}.
\]
The integral in the right hand side is easily seen to be $\preceq 2^{-kH/(H+1)}$. It follows that
\begin{equation}\label{firstterm}
\mathsf{E}\sup_{1\leq j\leq M_k}|\eta_k(\mathbf{x}_j)|\leq
2^{-kH/(2H+2)}\sqrt{1+\log M_k}.
\end{equation}

To estimate the second term we write
\[
\begin{aligned}
&\mathsf{E}\sup_{1\leq j\leq
M_k}\sup_{\mathbf{x},\mathbf{y}\in\mathcal{B}_j}|\eta_k(\mathbf{x})-\eta_k(\mathbf{y})|\\
&\quad\leq
\mathsf{E}\sup_{1\leq j\leq
M_k}\sup_{\mathbf{x},\mathbf{y}\in\mathcal{B}_j}\sum_{2^{k-1}<(m+1)(m/2+n)^{2H+1}\leq
2^k}\sum^{h(m,N)}_{l=1}|u^l_{mn}(\mathbf{x})-u^l_{mn}(\mathbf{y})|
\cdot|\xi^l_{mn}|.
\end{aligned}
\]
The difference $|u^l_{mn}(\mathbf{x})-u^l_{mn}(\mathbf{y})|$ can be estimated as
\[
\begin{aligned}
|u^l_{mn}(\mathbf{x})-u^l_{mn}(\mathbf{y})|&\leq|g_m(j_{|m-1|-H,n}\|\mathbf{x}\|)-\delta^0_m|
\cdot|S^l_m(\mathbf{x}/\|\mathbf{x}\|)-S^l_m(\mathbf{y}/\|\mathbf{y}\|)|\\
&\quad+|S^l_m(\mathbf{y}/\|\mathbf{y}\|)|\cdot|g_m(j_{|m-1|-H,n}\|\mathbf{x}\|)-
g_m(j_{|m-1|-H,n}\|\mathbf{y}\|)|.
\end{aligned}
\]
Let $\Delta_0$ denote the angular part of the Laplace operator in the space $\mathbb{R}^N$. Using the Mean Value Theorem for $g_m$, Lemma~\ref{l2} and formulas
\[
\begin{aligned}
\Delta_0S^l_m(\mathbf{x}/\|\mathbf{x}\|)&=-m(m+N-2)S^l_m(\mathbf{x}/\|\mathbf{x}\|),\\
|S^l_m(\mathbf{x}/\|\mathbf{x}\|)-S^l_m(\mathbf{y}/\|\mathbf{y}\|)|&\preceq
\sup_{\mathbf{z}}|\sqrt{-\Delta_0}S^l_m(\mathbf{z}/\|\mathbf{z}\|)|\varphi(\mathbf{x},\mathbf{y}),
\end{aligned}
\]
where the $\sup$ is taken over all points $\mathbf{z}$ belonging to the segment of the geodesic circle connecting the points $\mathbf{x}/\|\mathbf{x}\|$ and $\mathbf{y}/\|\mathbf{y}\|$, we obtain
\[
|u^l_{mn}(\mathbf{x})-u^l_{mn}(\mathbf{y})|\preceq\varepsilon_k\tau_{mn}\sqrt{m}(m/2+n).
\]
Hence, we have
\[
\mathsf{E}\sup_{1\leq j\leq
M_k}\sup_{\mathbf{x},\mathbf{y}\in\mathcal{B}_j}|\eta_k(\mathbf{x})-\eta_k(\mathbf{y})|\preceq
\varepsilon_k\sum_{2^{k-1}<(m+1)(m/2+n)^{2H+1}\leq
2^k}h(m,N)\tau_{mn}\sqrt{m}(m/2+n).
\]
Using Lemma~\ref{l2} once more, the last formula may be rewritten as
\[
\mathsf{E}\sup_{1\leq j\leq
M_k}\sup_{\mathbf{x},\mathbf{y}\in\mathcal{B}_j}|\eta_k(\mathbf{x})-\eta_k(\mathbf{y})|\preceq
\varepsilon_k\sum_{2^{k-1}<(m+1)(m/2+n)^{2H+1}\leq
2^k}m^{N-3/2}(m/2+n)^{1/2-H}.
\]
The integral, which corresponds to the sum in the right hand side, is easily seen to be asymptotically equal to $2^{kN/(2H+2)}$.
Therefore,
\[
\mathsf{E}\sup_{1\leq j\leq
M_k}\sup_{\mathbf{x},\mathbf{y}\in\mathcal{B}_j}|\eta_k(\mathbf{x})-\eta_k(\mathbf{y})|\preceq
\varepsilon_k\cdot 2^{kN/(2H+2)}.
\]
Combining this with the estimate \eqref{firstterm} for the first term, we obtain the asymptotic relation
\[
\mathsf{E}\sup_{\mathbf{x}\in\mathcal{B}}|\eta_k(\mathbf{x})|\preceq
2^{-kH/(2H+2)}\sqrt{1+\log M_k}+\varepsilon_k\cdot 2^{kN/(2H+2)}.
\]
Now set $\varepsilon_k=2^{-k(H+N)/(2H+2)}$ and recall that
$M_k\preceq\varepsilon^{-N}_k=2^{kN(H+N)/(2H+2)}$. Then we see that
the first term is bounded by a constant times
$2^{-kH/(2H+2)}\sqrt{k}$ and the second one is of lower order. We
proved that
\begin{equation}\label{mainestimate}
\mathsf{E}\sup_{\mathbf{x}\in\mathcal{B}}|\eta_k(\mathbf{x})|\preceq
2^{-kH/(2H+2)}\sqrt{k}.
\end{equation}

To complete the proof, it suffices to show that
\begin{equation}\label{last}
\mathsf{E}\sup_{\mathbf{x}\in\mathcal{B}}\left|\sum_{(m+1)(m/2+n)^{2H+1}>
q}\sum^{h(m,N)}_{l=1}u^l_{mn}(\mathbf{x})\xi^l_{mn}\right|\preceq
q^{-H/(2H+2)}(\log q)^{1/2}.
\end{equation}
Let $r$ be the positive integer such that $2^{r-1}<q\leq 2^r$. Then by the triangle inequality
\[
\begin{aligned}
&\left|\sum_{(m+1)(m/2+n)^{2H+1}>
q}\sum^{h(m,N)}_{l=1}u^l_{mn}(\mathbf{x})\xi^l_{mn}\right|\leq
\left|\sum_{(m+1)(m/2+n)^{2H+1}>
2^r}\sum^{h(m,N)}_{l=1}u^l_{mn}(\mathbf{x})\xi^l_{mn}\right|\\
&\quad+\left|\sum_{q<(m+1)(m/2+n)^{2H+1}\leq
2^r}\sum^{h(m,N)}_{l=1}u^l_{mn}(\mathbf{x})\xi^l_{mn}\right|\\
&\leq\sum_{k>r}|\eta_k(\mathbf{x})|+\left|\sum_{q<(m+1)(m/2+n)^{2H+1}\leq
2^r}\sum^{h(m,N)}_{l=1}u^l_{mn}(\mathbf{x})\xi^l_{mn}\right|.
\end{aligned}
\]
Using \eqref{mainestimate} and the fact that
\[
\sum_{k\geq r}\sqrt{k}2^{-kN/(2H+2)}\preceq\sqrt{r}2^{-rN/(2H+2)},
\]
we obtain
\[
\begin{aligned}
\sum_{k>r}\mathsf{E}\sup_{\mathbf{x}\in\mathcal{B}}|\eta_k(\mathbf{x})|&\preceq
\sum_{k>r}\sqrt{k}2^{-kN/(2H+2)}\\
&\preceq\sqrt{r-1}\cdot 2^{-(r-1)H/(2H+2)}\\
&\preceq q^{-H/(2H+2)}(\log q)^{1/2}.
\end{aligned}
\]
The arguments in the proof of \eqref{mainestimate} show that since $2^{r-1}<q\leq 2^r$, we also have
\[
\begin{aligned}
\mathsf{E}\sup_{\mathbf{x}\in\mathcal{B}}\left|\sum_{q<(m+1)(m/2+n)^{2H+1}\leq
2^r}\sum^{h(m,N)}_{l=1}u^l_{mn}(\mathbf{x})\xi^l_{mn}\right|&\preceq
\sqrt{r}2^{-rH/(2H+2)}\\
&\preceq q^{-H/(2H+2)}(\log q)^{1/2}.
\end{aligned}
\]
Relation \eqref{last} is proved.

\end{document}